\documentclass{article}
\usepackage{lineno,rotating,caption}
\usepackage{multirow}
\usepackage{authblk}
\usepackage{amsmath,amssymb,amsfonts,amsthm,anysize,setspace}
\usepackage{amssymb}
\usepackage{graphicx}
\usepackage[colorlinks]{hyperref}
\usepackage{caption}
\usepackage{graphics}
\usepackage{subfigure}
\usepackage{subfigure}
\usepackage{amsthm}
\usepackage{url}
\usepackage{pdftexcmds}
\usepackage{authblk}

\chardef\bslash=`\\ 

\newtheorem{thm}{Theorem}[section]
\theoremstyle{definition}

\newtheorem{rmk}{Remark}
\newtheorem{ex}{Example}

\begin{document}

\title{Modified least squares method and a review of its applications in machine learning and fractional differential/integral equations}

\author[1,2]{Abhishek Kumar Singh}
\affil[1]{Deptartment of Mathematics\\ Indian Institute of  Technology Delhi, India}
\affil[2]{Institute of Mathematics and Computer Science, Universit\"at Greifswald,
	Walther-Rathenau-Stra{\ss}e~47, 17489 Greifswald, Germany}

\author[1]{Mani Mehra\footnote{Corresponding author}}
\author[3]{Anatoly A. Alikhanov}
\affil[3]{North-Caucasus Center for Mathematical Research, North-Caucasus Federal University, Russia \protect\\
\texttt{assinghabhi@gmail.com, mmehra@maths.iitd.ac.in, aaalikhanov@gmail.com}}

\maketitle




\begin{abstract}
The least squares method provides the best-fit curve by minimizing the total squares error. In this work, we provide the modified least squares method based on the fractional orthogonal polynomials that belong to the space $M_{n}^{\lambda} := \text{span}\{1,x^{\lambda},x^{2\lambda},\ldots,x^{n\lambda}\},~\lambda \in (0,2]$. Numerical experiments demonstrate how to solve different problems using the modified least squares method. Moreover, the results show the advantage of the modified least squares method compared to the classical least squares method. Furthermore, we discuss the various applications of the modified least squares method in the fields like fractional differential/integral equations and machine learning.
\end{abstract}

\numberwithin{equation}{section}
\section{Introduction}\label{6s1}
The least squares method is one of the oldest methods of modern statistics used to obtain the physical parameters from the experimental data. The first use of the least squares method is generally attributed to Gauss in 1795, although Legendre concurrently and independently used it \cite{harter1974method}. Gauss invented the least squares method to estimate planets' orbital motion from telescopic measurements. In modern statistics, Galton \cite{abdi2007method} was the first to use the least squares method in his work on the heritability of size, which laid down the foundations of correlation and regression analysis. 
Nowadays, the least squares method is widely used to find the best-fit curve while finding the parameter involved in the curve. There are many versions of the least squares method available in the literature. The simpler version is called the ordinary least squares method, and the more advanced one is the weighted least squares method, which performs better than the ordinary least squares method. The recent version of the least squares method is the moving least squares method \cite{bale2021one}, and the partial least squares method \cite{smoliak2010application}.\par 
One of the areas where the least squares method is frequently used is machine learning, where we analyze data for regression analysis and classification \cite{raissi2018hidden}. Machine learning is a field of artificial intelligence that allows computer systems to learn using available data. Recently machine learning algorithms (regression analysis and classification) have become very popular for analyzing data and making predictions. Another application of the least squares method is solving fractional differential/integral equations. Fractional differential/integral equations give an excellent way to deal with complex phenomena in nature, such as biological systems, control theory, finance, signal and image processing, sub-diffusion and super-diffusion process, viscoelastic fluid, electrochemical processes, and so on \cite{alikhanov2017time,herrmann2013folded,srivastava2021efficient,singh2021wavelet,mehandiratta2019existence}. The fractional differential equations are equivalent to the Hammerstein form of Volterra's second kind integral equations for the specific choice of kernel (for more details see \cite{diethelm2012volterra}). Due to the importance of fractional differential/integral equations, people are interested in solving them numerically because of the non-availability of exact solutions. Many numerical methods are available in the literature to solve fractional differential/integral equations, such as finite difference \cite{meerschaert2006finite,mehandiratta2021optimal}, compact finite difference \cite{patel2020fourth,alikhanov2021crank}, finite element \cite{deng2009finite,jin2015galerkin} , and wavelet methods \cite{kumar2021legendre,singh2020uncertainty}. Recently, the least squares method based on classical polynomials is also used to solve fractional differential/integral equations \cite{rahimkhani2017fractional}. \par 
The aforementioned discussion concludes with the importance of the least squares method in various disciplines. In this work, we proposed the modified least squares method based on fractional polynomials. The main contribution of this paper consists of the following aspects:
\begin{itemize}
	\item A modified least squares method based on fractional orthogonal polynomials belongs to the space $M_{n}^{\lambda} := \text{span}\{1,x^{\lambda},x^{2\lambda},\ldots,x^{n\lambda}\},~\lambda \in (0,2]$  has been proposed.
	\item Applications of the modified least squares method have been discussed in detail, especially in the field of fractional differential/integral  equations and machine learning.
\end{itemize}
The paper has been arranged in the following pattern: Section \ref{6s2} describes some basic concepts of approximation theory and orthogonal polynomials like M\"untz-Legendre polynomials. Section \ref{6s3} deals with the review of a least squares method. In Section \ref{6s4}, we develop a modified least squares method based on fractional orthogonal polynomials. Section \ref{6s5} provides the numerical results to validate the modified least squares method. Section \ref{6s6} summarizes the applications part of the method. Finally, Section \ref{6s7} gives the brief conclusion.

\section{\textbf{Preliminaries}}\label{6s2}
In this Section, some necessary results from approximation theory and orthogonal polynomials are summarized. Let define the space
\begin{equation}\label{61e1}
M^{\lambda} := \text{span}\{x^{\lambda_{0}},x^{\lambda_{1}},\ldots,x^{{\lambda_{n}}},\ldots\},
\end{equation}
and for its subspaces we define
\begin{equation}\label{61e2}
M_{n}^{\lambda} := \text{span}\{x^{\lambda_{0}},x^{\lambda_{1}},\ldots,x^{{\lambda_{n}}}\},
\end{equation}
with $0 \leq \lambda_{0} < \lambda_{1} < \ldots \rightarrow \infty.$ The space $M^{\lambda}$ is known as M\"untz space. Now, we recall one of the fundamental theorems of approximation theory called M\"untz-Sz$\acute{a}$sz theorem, which is related to the denseness of polynomial belonging to space $M^{\lambda}$.
\begin{thm}{(M\"untz-Sz$\acute{a}$sz Theorem)}
The M\"untz polynomials of the forms $\sum_{i=0}^{n}a_{i}x^{\lambda_{i}} \in M_{n}^{\lambda}$ with real coefficients are dense in $L^{2}[0,1]$ if and only if
\begin{equation}\label{61e3}
\sum_{i=1}^{\infty}\frac{1}{\lambda_{i}} = +\infty.
\end{equation}	
Moreover, if $\lambda_{0} = 0$, then M\"untz polynomials are dense in $C[0,1]$.
\end{thm}
\noindent \textbf{Proof.} See \cite{borwein1994muntz}.\par
In this work, we assume that $\lambda_{i} = i\lambda,~i=0,1,\ldots,n,\ldots$, where $\lambda$ is a real constant and $\lambda \in (0,2]$. Now, we review one of the general forms of classical orthogonal polynomials called Jacobi polynomials. Moreover, we also review the  M\"untz-Legendre polynomials, and because of some advantage in terms of computational accuracy, we define the M\"untz-Legendre polynomials with the help of Jacobi polynomials.
\subsection{Jacobi polynomials}
The Jacobi polynomials with parameters $a,b > -1$, denoted by $P_{n}^{a,b}$, is defined in the interval $-1 \leq x \leq 1$ as
\begin{equation}\label{42e1}
P_{n}^{a,b}(x) = \sum_{i = 0}^{n}\frac{(-1)^{n-i}(1+b)_{n}(1+a+b)_{i+n}}{i!(n-i)!(1+b)_{i}(1+a+b)_{n}}\left(\frac{1+x}{2}\right)^{i},
\end{equation}
where $(1+b)_{i} = (1+b)(2+b)\ldots(i+b)$ and $(1+b)_{0} = 1$. \par
In practice, one can compute the Jacobi polynomials using the following recurrence relation 
\begin{equation}\label{42e2}
\begin{aligned}
P_{0}^{a,b}(x)  = 1,&~~ P_{1}^{a,b}(x) = \frac{1}{2}[(a-b)+(a+b+2)x], \\
c_{1,n}^{a,b}P_{n}^{a,b}(x) & = c_{2,n}^{a,b}(x)P_{n}^{a,b}(x)-c_{3,n}^{a,b}P_{n-1}^{a,b}(x),
\end{aligned}
\end{equation}
where
\begin{equation}\label{42e3}
\begin{aligned}
c_{1,n}^{a,b}& = 2(n+1)(n+a+b+1)(2n+a+b), \\
c_{2,n}^{a,b}(x) & = (2n+a+b+1)[(2n+a+b)(2n+a+b+2)x+a^{2}-b^{2}], \\
c_{3,n}^{a,b} & = 2(n+a)(n+b)(2n+a+b+2).
\end{aligned}
\end{equation}
\subsection{M\"untz Legendre polynomials}
One can defined M\"untz Legendre polynomials on the interval $[0,1]$ as follows
\begin{equation}\label{42e4}
L_{n}(x) = \sum_{i = 0}^{n}\eta_{n,i}x^{\lambda_{i}},~~\eta_{n,i} = \frac{\prod_{k = 0}^{n-1}(\lambda_{i}+\bar{\lambda}_{k}+1)}{\prod_{k = 0,k \neq i}^{n}(\lambda_{i}-\lambda_{k})}.
\end{equation}
These polynomials satisfy the following orthogonality condition with respect to weight function $W(x) = 1$
\begin{equation}
\int_{0}^{1}L_{n}(x)L_{m}(x)dx = \frac{\delta_{n,m}}{\lambda_{n}+\bar{\lambda}_{n}+1}.
\end{equation}
Since we assumed that $\lambda_{i} = i\lambda$, where $\lambda$ is a real constant, then the M\"untz-Legendre polynomials on the interval $[0, 1]$ are represented by the formula
\begin{equation}
L_{n}(x ; \lambda) := \sum_{i = 0}^{n} \eta_{n,i} x^{i\lambda},~~ \eta_{n,i} = \frac{(-1)^{n-i}}{\lambda^{n}i!(n-i)!}\prod_{k = 0}^{n-1}((i+k)\lambda+1).
\end{equation}
From Equation \eqref{42e4}, one can observe that evaluating M\"untz-Legendre polynomials, mainly when $n$ is vast, and $x$ is closed to $1$, is problematic in finite arithmetic. Milovanovic has addressed these problems \cite{milovanovic1999muntz}. We will use a method for evaluating M\"untz-Legendre polynomials, which is based on three-term recurrence relation induced from the accompanying theorem with the help of Jacobi polynomials.
\begin{thm}\label{6thm1}
	Let $\lambda > 0$ be a real number and $x \in [0, 1]$. Then the representation 
	$$L_{n}(x; \lambda) = P_{n}^{(0,1/\lambda-1)}\left(2x^{\lambda}-1\right)$$
	holds true.
\end{thm}
\noindent \textbf{Proof.} See \cite{esmaeili2011numerical}. \hfill $\square$ \par
So in view of Theorem \ref{6thm1} and Equations \eqref{42e2} and \eqref{42e4}, the M\"untz-Legendre polynomials $L_{n}(x; \lambda)$ can be evaluated by means of three-terms recurrence relation 
\begin{equation}\label{42e5}
\begin{aligned}
L_{0}(x; \lambda) = 1,~~ L_{1}(x; \lambda) = \left(\frac{1}{\lambda}+1\right)x^{\lambda}-\frac{1}{\lambda}, \\
d_{1,n}L_{n+1} = d_{2,n}(t)L_{n}(x; \lambda)-d_{3,n}L_{n-1}(x; \lambda),
\end{aligned}
\end{equation}
where
$$d_{1,n} = c_{1,n}^{0,1/\lambda-1},~~ d_{2,n}(x) = c_{2,n}^{0,1/\lambda-1}\left(2x^{\lambda}-1\right),~~d_{3,n} = c_{3,n}^{0,1/\lambda-1}.$$ 
\section{Least squares method}\label{6s3}
For approximating the continuous function $y(x)$ defined on the interval $[a,b]$ with an algebraic polynomial $P_{n}(x;1) = \sum_{i = 0}^{n}a_{i}x^{i} \in M_{n}^{1}$ using least squares method, we choose the constants $a_{0},a_{1},\dots,a_{n}$, which minimize the least squares error $E^{C}(a_{0},a_{1},\dots,a_{n})$, where
\begin{equation*}
E^{C}(a_{0},a_{1},\dots,a_{n}) = \int_{a}^{b}(y(x)-P_{n}(x))^{2}dx.
\end{equation*}
A necessary condition for real coefficient $a_{0},a_{1},\ldots,a_{n}$ that minimize the error $E^{C}(a_{0},a_{1},\dots,a_{n})$ is that
\begin{equation*}
\frac{\partial E^{C}}{\partial a_{i}} = 0,~\text{for each}~i=0,1,\ldots,n.
\end{equation*}
Thus, one can get the following linear system of equations
\begin{equation}\label{62e1}
A_{LSM}^{C}X = b^{C},
\end{equation}
where $X = [a_{0},a_{1},\ldots,a_{n}]^{T},~A_{LSM}^{C}$ is the $(n+1)\times(n+1)$ coefficient matrix and $b^{C}$ is the $(n+1)\times 1$ column vector. The $(i,j)^{th}$ element of the matrix $A_{LSM}^{C}$ and $i^{th}$ element of the column vector $b^{C}$ are given by
$$(A_{LSM}^{C})_{i,j} = \int_{a}^{b}x^{i+j-2}dx,~~\text{and},~~ (b^{C})_{i} = \int_{a}^{b}y(x)x^{i-1}dx,$$
respectively. Similarly, in discrete case, for approximating a data set $(x_{k},y_{k}),~k = 1,2,\dots,m$, with an algebraic polynomial $P_{n}(x;1) \in M_{n}^{1}$ using least squares method, we choose the constants $a_{0},a_{1},\dots,a_{n}$ which minimize the least squares error $E^{D}(a_{0},a_{1},\dots,a_{n})$, where
\begin{equation*}
\begin{aligned}
E^{D}(a_{0},a_{1},\dots,a_{n}) &= \sum_{k=1}^{m}(y_{k}-P_{n}(x_{k}))^{2}, \\
& = \sum_{k = 1}^{m}y_{k}^{2}-2\sum_{i=0}^{n}a_{i}(\sum_{k=1}^{m}y_{k}x_{k}^{i})+2\sum_{i=0}^{n}\sum_{j=0}^{n}a_{i}a_{j}(\sum_{k=1}^{m}x_{k}^{i+j}).
\end{aligned}
\end{equation*}
For $E^{D}$ to be minimized it is necessary that $\frac{\partial E^{D}}{\partial a_{i}} = 0$, for each $i=0,1,\dots,n.$ Thus, one can get following linear system of equations
\begin{equation}\label{62e2}
A_{LSM}^{D}X = b^{D},
\end{equation}
where $A_{LSM}^{D}$ is the $(n+1)\times(n+1)$ coefficient matrix and $b^{D}$ is the $(n+1)\times 1$ column vector. The $(i,j)^{th}$ element of the matrix $A_{LSM}^{D}$ and $i^{th}$ element of the column vector $b^{D}$ are given by
$$(A_{LSM}^{D})_{i,j} = \sum_{k=1}^{m}x_{k}^{i+j-2},~~\text{and},~~ (b^{D})_{i} = \sum_{k=1}^{m}y_{k}x_{k}^{i-1},$$
respectively.
\section{Modified least squares method}\label{6s4}
In this Section, we proposed a modified least squares method based on a fractional polynomials that belongs to the space $M_{n}^{\lambda}$.
Apply a similar technique which is discussed in Section \ref{6s3} for any algebraic polynomial $P_{n}(x;\lambda) = \sum_{i = 0}^{n}a_{i}x^{i\lambda} \in M_{n}^{\lambda}$ to approximating the continuous function $y(x)$ defined on the interval $[a,b]$, we choose the constants $a_{0},a_{1},\dots,a_{n}$ which minimize the least squares error $E^{C}(a_{0},a_{1},\dots,a_{n};\lambda)$, where
\begin{equation}\label{newe1}
E^{C}(a_{0},a_{1},\dots,a_{n};\lambda) = \int_{a}^{b}(y(x)-P_{n}(x;\lambda))^{2}dx.
\end{equation}
A necessary condition for real coefficient $a_{0},a_{1},\ldots,a_{n},$ that minimize the error $E^{C}(a_{0},a_{1},\dots,a_{n};\lambda)$ is that
\begin{equation*}
\frac{\partial E^{C}}{\partial a_{i}} = 0,~\text{for each}~i=0,1,\ldots,n.
\end{equation*}
Thus, one can get the following linear system of equations
\begin{equation}\label{63e1}
A_{MLSM}^{C}X = d^{C},
\end{equation}
where $A_{MLSM}^{C}$ is the $(n+1)\times(n+1)$ coefficient matrix and $d^{C}$ is the $(n+1)\times 1$ column vector. The $(i,j)^{th}$ element of the matrix $A_{MLSM}^{C}$ and $i^{th}$ element of the column vector $d^{C}$ are given by
$$(A_{MLSM}^{C})_{i,j} = \int_{a}^{b}x^{\lambda_{i-1}+\lambda_{j-1}}dx,~~\text{and},~~ (d^{C})_{i} = \int_{a}^{b}y(x)x^{\lambda_{i-1}}dx,$$
respectively.
Similarly, in discrete case, for approximating the data set $(x_{k},y_{k}),~k = 1,2,\ldots,m,$ with an algebraic polynomial $P_{n}(x;\lambda)$ using least squares method, we choose the constants $a_{0},a_{1},\dots,a_{n}$, which minimize the least square error $E^{D}(a_{0},a_{1},\dots,a_{n};\lambda)$, where
\begin{equation}\label{newe2}
\begin{aligned}
E^{D}(a_{0},a_{1},\dots,a_{n};\lambda) &= \sum_{k = 1}^{m}(y_{k}-P_{n}(x_{k};\lambda))^{2}, \\
& = \sum_{k=1}^{m}y_{k}^{2}-2\sum_{k=1}^{m}P_{n}(x_{k};\lambda)y_{k}+\sum_{k=1}^{m}(P_{n}(x_{k};\lambda))^{2},\\
& = \sum_{k = 1}^{m}y_{k}^{2}-2\sum_{i = 0}^{n}a_{i}(\sum_{k=1}^{m}y_{k}x_{k}^{\lambda_{i}})+2\sum_{i =0}^{n}\sum_{j = 0}^{n}a_{i}a_{j}(\sum_{k = 1}^{m}x_{k}^{\lambda_{i}+\lambda_{j}}).
\end{aligned}
\end{equation}
For $E^{D}$ to be minimized it is necessary that $\frac{\partial E^{D}}{\partial a_{i}} = 0$, for each $i=0,1,\dots,n.$ Thus, one can get following linear system of equations
\begin{equation}\label{63e2}
A_{MLSM}^{D}X = d^{D},
\end{equation}
where $A_{MLSM}^{D}$ is the $(n+1)\times(n+1)$ coefficient matrix and $d^{D}$ is the $(n+1)\times 1$ column vector. The $(i,j)^{th}$ element of the matrix $A_{MLSM}^{D}$ and $i^{th}$ element of the column vector $d^{D}$ are given by
$$(A_{MLSM}^{D})_{i,j} = \sum_{k=1}^{m}x_{k}^{\lambda_{i-1}+\lambda_{j-1}},~~\text{and},~~ (d^{D})_{i} = \sum_{k=1}^{m}y_{k}x_{k}^{\lambda_{i-1}},$$
respectively.
We know that the large value of $n$, the matrix $A_{MLSM}^{C}$ and $A_{MLSM}^{D}$ become ill-conditioned, which causes significant errors in estimating the parameters $a_{i},i=0,1,\ldots,n$. This difficulty can be avoided if the functions belonging to the space $M_{n}^{\lambda}$, denoted by $L_{i}(x;\lambda,W),~i=0,1,\ldots,n$, are so chosen that they are orthogonal with respect to the weight function $W(x)$ over the interval $[a,b]$. In this case the error function in the continuous and discrete case becomes
\begin{equation}\label{63e3}
E^{C}(a_{0},a_{1},\dots,a_{n};\lambda) = \int_{a}^{b}W(x)(y(x)-\sum_{i=0}^{n}a_{i}L_{i}(x;\lambda,W))^{2}dx,
\end{equation}
and 
\begin{equation}\label{63e4}
E^{D}(a_{0},a_{1},\dots,a_{n};\lambda) = \sum_{k = 1}^{m}W(x_{k})(y_{k}-\sum_{i=0}^{n}a_{i}L_{i}(x_{k};\lambda,W))^{2},
\end{equation}
respectively. 
The necessary condition for real coefficients $a_{0},a_{1},\ldots,a_{n}$ which minimize the error gives the normal equations. The normal equations in the continuous and discrete cases are
\begin{equation}\label{63e7}
\int_{a}^{b}W(x)(y(x)-\sum_{i=0}^{n}a_{i}L_{i}(x;\lambda,W))L_{j}(x;\lambda,W)dx = 0,~j=0,1,\ldots,n,
\end{equation}
and 
\begin{equation}\label{63e8}
\sum_{k=1}^{m}W(x_{k})(y_{k}-\sum_{i=0}^{n}a_{i}L_{i}(x_{k};\lambda,W))L_{j}(x_{k};\lambda,W) = 0,~j=0,1,\ldots,n,
\end{equation}
respectively. One can observe from the Equations \eqref{63e7} and \eqref{63e8} the parameter value can be determined directly. Thus the use of orthogonal functions not only avoids the problem of ill-conditioning but also determines the constants $a_{i},~i=0,1,\ldots,n,$ directly. As discussed in Section \ref{6s2}, the M\"untz-Legendre polynomials $L_{i}(x;\lambda) \in M_{n}^{\lambda},~i=0,1,\ldots,n,$ are orthogonal for the weight function $W(x) = 1$ on the interval $[0,1]$. If we consider the following error function in the continuous and discrete case
\begin{equation}\label{63e5}
E^{C}(a_{0},a_{1},\dots,a_{n};\lambda) = \int_{0}^{1}W(x)(y(x)-\sum_{i=0}^{n}a_{i}L_{i}(x;\lambda))^{2}dx,
\end{equation}
and 
\begin{equation}\label{63e6}
E^{D}(a_{0},a_{1},\dots,a_{n};\lambda) = \sum_{k = 1}^{m}W(x_{k})(y_{k}-\sum_{i=0}^{n}a_{i}L_{i}(x_{k};\lambda))^{2},
\end{equation}
respectively. One can easily observe that if we find the necessary condition for the parameters $a_{0},a_{1},\ldots,a_{n},$ we get the systems of equations. For finding the value of $a_{i},~i=0,1,\ldots,n,$ we have to solve the equation system because M\"untz-Legendre polynomials are not orthogonal with respect to the weight function $W(x) \neq 1$. Therefore, one may face the problem of ill-conditioning for a large value of $n$. This problem can be avoided if we generate the orthogonal polynomials with respect to the weight function $W(x)$.
\begin{rmk}
The modified least squares method is based on the fractional polynomials, which depend on the fractional parameter $\lambda$. If we put $\lambda = 1$ in the Equations \eqref{63e1} and \eqref{63e2} then we get \eqref{62e1} and \eqref{62e2} respectively. Therefore, in the modified least squares method, we have one additional degree of freedom compared to the classical least squares method in the form of fractional parameter $\lambda$. So, we tune the parameter $\lambda$ to get the desired results.
\end{rmk}	
\begin{rmk}\label{63rmk1}
	Since M\"untz-Legendre polynomials are orthogonal with respect to weight function $W(x) = 1$. We know that the fractional derivative/integral of the polynomial belonging to space $M_{n}^{\lambda}$ is again in space $M_{n}^{\lambda}$ when $\lambda$ and fractional order are the same. Hence, while solving the fractional differential equations using the least squares method and avoiding the difficulty due to ill-conditioning, we need polynomials belonging to the space $M_{n}^{\lambda}$, which are orthogonal to weight function $(b-x)^{-\lambda}$ or $(x-a)^{-\lambda}$. Moreover, usually, in applications, only a part of the given data needs more attention; for example, in some cases, the data may have more accuracy in some regions than in others. In such cases, the weight function indicates where data should be given more importance, and it should be chosen accordingly.
\end{rmk}
Now, we are going to discuss the Theorem, which helps us to generate the fractional orthogonal polynomials belonging to the space $M_{n}^{\lambda}$ with respect to the weight function $W(x)$.
\begin{thm}\label{64thm1}
The set of fractional polynomial functions $\{L_{0}(x;\lambda,W),L_{1}(x;\lambda,W),\ldots,L_{n}(x;\lambda,W)\}$ defined in the following way is orthogonal on $[a,b]$ with respect to the weight function $W(x)$. 
\begin{equation*}
		L_{0}(x;\lambda,W) = 1,~~ L_{1}(x;\lambda,W) = x^{\lambda}-B_{1},
\end{equation*}	
where 
\begin{equation*}
	B_{1} = \frac{\int_{a}^{b}W(x)x^{\lambda}(L_{0}(x;\lambda,W))^{2}dx}{\int_{a}^{b}W(x)(L_{0}(x;\lambda,W))^{2}dx},
\end{equation*}
and when $i \geq 2$,
\begin{equation*}
	L_{i}(x;\lambda,W) = (x^{\lambda}-B_{i})L_{i-1}(x;\lambda,W)-C_{i}L_{i-2}(x;\lambda,W),
\end{equation*}
where 
\begin{equation*}
	\begin{aligned}
		B_{i} &= \frac{\int_{a}^{b}W(x)x^{\lambda}(L_{i-1}(x;\lambda,W))^{2}dx}{\int_{a}^{b}W(x)(L_{i-1}(x;\lambda,W))^{2}dx},~~ C_{i} = \frac{\int_{a}^{b}W(x)x^{\lambda}L_{i-1}(x;\lambda,W)L_{i-2}(x;\lambda,W)dx}{\int_{a}^{b}W(x)(L_{i-2}(x;\lambda,W))^{2}dx}.
	\end{aligned}
\end{equation*}
\end{thm}
\noindent \textbf{Proof.} We prove that the above result holds by mathematical induction; firstly, we will consider $L_{0}(x;\lambda,W)$ and $L_{1}(x;\lambda,W)$ and show that they are orthogonal.
\begin{equation*}
	\begin{aligned}
		\int_{a}^{b}W(x)L_{0}(x;\lambda,W)L_{1}(x;\lambda,W)dx & = \int_{a}^{b}W(x)(x^{\lambda}-B_{1})dx \\
		&= \int_{a}^{b}W(x)x^{\lambda}dx-B_{1}\int_{a}^{b}W(x)dx \\
		& = 0.
	\end{aligned}
\end{equation*}
Now, For the induction hypothesis, assume that $\{L_{0}(x;\lambda,W),L_{1}(x;\lambda,W),\ldots,L_{n-1}(x;\lambda,W)\}$ is orthogonal on $[a,b]$ with respect to the weight function $W(x)$. Consider the following
\begin{equation*}
	\int_{a}^{b}W(x)L_{i}(x;\lambda,W)L_{n}(x;\lambda,W)dx,~i=0,1,\ldots,n-1.
\end{equation*}
If we will prove the above integral value is zero then we are done. Take $i = n-1$ in the above expression and using the recurrence relation for $L_{n}(x;\lambda,W)$, we get
\begin{equation*}
	\begin{aligned}
		\int_{a}^{b}W(x)L_{n-1}(x;\lambda,W)L_{n}(x;\lambda,W)dx & = \int_{a}^{b}W(x)x^{\lambda}L_{n-1}^{2}(x;\lambda,W)dx-B_{n}\int_{a}^{b}W(x)L_{n-1}^{2}(x;\lambda,W)dx \\
		&= 0.
	\end{aligned}
\end{equation*}
Similarly, we can show that integral value is zero for $i = n-2$. Now, consider the integral for $i = n-3$ and using the recurrence relation for $L_{n}(x;\lambda,W)$, we get
\begin{equation}\label{6eth1}
	\begin{aligned}
		\int_{a}^{b}W(x)L_{n-3}(x;\lambda,W)L_{n}(x;\lambda,W)dx &= \int_{a}^{b}W(x)x^{\lambda}L_{n-3}(x;\lambda,W)L_{n-1}(x;\lambda,W)dx.
	\end{aligned}
\end{equation}
Using the relation $x^{\lambda}L_{n-3}(x;\lambda,W) = L_{n-2}(x;\lambda,W)+B_{n-2}L_{n-3}(x;\lambda,W)+C_{n-2}L_{n-4}(x;\lambda,W)$ and orthogonal property in the Equation \eqref{6eth1}, we get
\begin{equation*}
	\int_{a}^{b}W(x)L_{n-3}(x;\lambda,W)L_{n}(x;\lambda,W)dx = 0.
\end{equation*}
From the similar argument, we can show for $i=0,1,2,\ldots,n-4$. Therefore, the set of fractional polynomial functions $\{L_{0}(x;\lambda,W),L_{1}(x;\lambda,W),\ldots,L_{n}(x;\lambda,W)\}$ defined in the following way is orthogonal on $[a,b]$ with respect to the weight function $W(x)$. \hfill $\square$ \par We can generalize the above Theorem in the discrete case. The following Theorem will help us to generate the fractional orthogonal polynomials over a set of points $x_{k},~k = 1,2,\ldots,m$ with respect to the weight function $W(x)$.
\begin{thm}\label{64thm2}
The set of fractional polynomial functions $\{L_{0}(x;\lambda,W),L_{1}(x;\lambda,W),\ldots,L_{n}(x;\lambda,W)\}$ defined in the following way is orthogonal over a set of points $x_{k},~k = 1,2,\ldots,m$, with respect to the weight function $W(x)$.
\begin{equation*}
	L_{0}(x;\lambda,W) = 1,~~ L_{1}(x;\lambda,W) = x^{\lambda}-B_{1},
\end{equation*}	
where 
\begin{equation*}
	B_{1} = \frac{\sum_{k=1}^{m}W(x_{k})x_{k}^{\lambda}(L_{0}(x_{k};\lambda,W))^{2}}{\sum_{k=1}^{m}W(x_{k})(L_{0}(x_{k};\lambda,W))^{2}},
\end{equation*}
and when $i \geq 2$,
\begin{equation*}
	L_{i}(x;\lambda,W) = (x^{\lambda}-B_{i})L_{i-1}(x;\lambda,W)-C_{i}L_{i-2}(x;\lambda,W),
\end{equation*}
where 
\begin{equation*}
	\begin{aligned}
		B_{i} &= \frac{\sum_{k=1}^{m}W(x_{k})x_{k}^{\lambda}(L_{i-1}(x_{k};\lambda,W))^{2}}{\sum_{k=1}^{m}W(x_{k})(L_{i-1}(x_{k};\lambda,W))^{2}},~~ C_{i} = \frac{\sum_{k=1}^{m}W(x_{k})x_{k}^{\lambda}L_{i-1}(x_{k};\lambda,W)L_{i-2}(x_{k};\lambda,W)}{\sum_{k=1}^{m}W(x_{k})(L_{i-2}(x_{k};\lambda,W))^{2}}.
	\end{aligned}
\end{equation*}	
\end{thm}
\noindent \textbf{Proof.} Proof is similar to the proof of the Theorem \ref{64thm1}, but we take discrete inner product in this case.\hfill $\square$ \par
After generating the fractional orthogonal polynomials with respect to weight function $W(x)$ using the Theorems \ref{64thm1} and \ref{64thm2}, we can directly compute the value of $a_{i}$. Hence, while minimizing the error defined by the Equations \eqref{63e3} and \eqref{63e4}, we get
\begin{equation}\label{63e9}
	a_{i} = \frac{\int_{a}^{b}W(x)y(x)L_{i}(x;\lambda,W)dx}{\int_{a}^{b}W(x)(L_{i}(x;\lambda,W))^{2}dx},~~i = 0,1,\ldots,n,
\end{equation} 
and 
\begin{equation}\label{63e10}
	a_{i} = \frac{\sum_{k=1}^{m}W(x_{k})y(x_{k})L_{i}(x_{k};\lambda,W)}{\sum_{k=1}^{m}W(x_{k})(L_{i}(x_{k};\lambda,W))^{2}},~~i = 0,1,\ldots,n,
\end{equation}
respectively.
\section{Test examples}\label{6s5}
This Section is devoted to illustrate the accuracy and efficiency of the proposed modified least squares method discussed in Section \ref{6s4}. All the numerical simulation were run on an Intel Core $i5-1135G7$, $2.42GHz$ machine with $16$GB RAM. To demonstrate the efficiency of the proposed method, we consider the following Examples:
\begin{ex}\label{6ex1}
	In this Example, we consider the function $y(x) = x^{0.75}+x^{1.5}$ defined on the interval $[0,1]$.
\end{ex}
Here, we approximate $y(x)$ with the $P_{2}(x;\lambda) = a_{0}+a_{1}x^{\lambda}+a_{2}x^{2\lambda} \in M_{2}^{\lambda}$. After implementation of the proposed method discussed in Section \ref{6s4}, the results are described in detail as below:
\begin{itemize}
	\item One can clearly observed that, for choises of $a_{0} = 0,~a_{1} = 1$ and $a_{2} = 1$ with $\lambda = 0.75$, we get the $P_{2}(x;0.75)$ exact as a $y(x)$.
	\item Table \ref{6tb1} shows the least squares error $E^{C}$ by the proposed method and CPU time with different value of $\lambda$.
	\item Also, Table \ref{6tb1} shows our proposed method's well accurate than the classical least squares method due to the additional parameter $\lambda$.
\end{itemize}
\begin{table}[ht]
	\begin{center}
		\begin{tabular}{| c| l| c|c|}
			\hline
			 $\lambda$ & $\boldsymbol{y(x) = x^{0.75}+x^{1.5}}$ & Total Error $E^{C}$ & CPU time (in second)  \\
			 \hline
			 0.75 & $y(x) \approx 0.0000+1.0000x^{0.75}+1.0000x^{1.5}$ & 2.70$e$-24 & 0.0558 \\
			\hline 
			 1 & $y(x) \approx 0.0329+1.7039x+0.2597x^{2}$ & 1.40$e$-5 & 0.0702 \\
			\hline 
			 1.5 & $y(x) \approx 0.1388+2.5269x^{1.5}-0.7126x^{3}$ & 8.78$e$-4 & 0.0631\\
			\hline
		\end{tabular}
		\caption{Comparison of the least squares error $E^{C}$ in Example \ref{6ex1} corresponding to $n = 2$ and  with different value of $\lambda$.}\label{6tb1}
	\end{center}
\end{table}
\begin{ex}\label{6ex6}
	In this Example, we consider the function $y(x) = x^{0.5}-\frac{\pi}{4}$ defined on the interval $[0,1]$.
\end{ex}
Here, we consider the two cases : first we approximate $y(x)$ with the  $a_{0}L_{0}(x;\lambda,W)+a_{1}L_{1}(x;\lambda,W)$, where $W(x) = (1-x)^{-\lambda}$, while in second case we approximate $y(x)$ with $b_{0}L_{0}(x;\lambda)+b_{1}L_{1}(x;\lambda)$. Both the cases, we will consider the followings residual error
\begin{equation*}
	\begin{aligned}
		E^{C}(a_{0},a_{1};\lambda) &= \int_{0}^{1}W(x)(y(x)-\sum_{i = 0}^{1}a_{i}L_{i}(x;\lambda,W))^{2}dx,
	\end{aligned}
\end{equation*}
and 
\begin{equation*}
	\begin{aligned}
		E^{C}(b_{0},b_{1};\lambda) &= \int_{0}^{1}W(x)(y(x)-\sum_{i = 0}^{1}b_{i}L_{i}(x;\lambda))^{2}dx,
	\end{aligned}
\end{equation*}
respectively. After implementation of the proposed method discussed in Section \ref{6s4}, the results are described in detail as below:
\begin{itemize}
	\item For $\lambda = 0.5$, if one can generate the fractional orthogonal polynomials up to order $1$ with respect to weight function $W(x) = (1-x)^{-0.5}$ using Theorem \ref{64thm1}, they get $L_{0}(x;0.5,W) = 1$ and $L_{1}(x;0.5,W) = x^{0.5}-\frac{\pi}{4}$.
	\item In first case, one can clearly observed that, for choices of $a_{0} = 0$ and $a_{1} = 1$ with $\lambda = 0.5$, we get the $a_{0}L_{0}(x;0.5,W)+a_{1}L_{1}(x;0.5,W)$ exact as a $y(x)$, without solving system of equations.
	\item In second cases, we get the value of $b_{0} = -0.1187$ and $b_{1} = 0.3333$, with solving the system of equations and $b_{0}L_{0}(x;0.5)+b_{1}L_{1}(x;0.5)$ exact as a $y(x)$. 
	\item From the discussion of the above results, we can conclude that generating the orthogonal polynomial with respect to weight function $W(x)$ gives the result without solving the system of equations and avoiding the ill-conditioning situation.
\end{itemize}
\begin{ex}\label{6ex2}
	In this Example, we consider the function $y(x) = x^{1.5}$ defined on the interval $[10,20]$ to generate the data set of length $20$.
\end{ex}
Here, we approximate $y(x)$ with the $P_{1}(x;\lambda) = a_{0}+a_{1}x^{\lambda} \in M_{1}^{\lambda}$. After implementation of the proposed method discussed in Section \ref{6s4}, the results are described in detail as below:
\begin{itemize}
	\item One can clearly observed that, for choices of $a_{0} = 0$ and $a_{1} = 1$ with $\lambda = 1.5$, we get the $P_{1}(x;1.5)$ exact as a $y(x)$.
	\item Table \ref{6tb2} shows the least squares error $E^{D}$ by the proposed method with different value of $\lambda$ and CPU time taken by the proposed method.
	\item Also, Table \ref{6tb2} displays the results for the noisy data set with generated by the function $y(x)$ with $5\%$ and $10\%$ noise. Table results show the our proposed method is robust to noise.
	\item From Table \ref{6tb2}, one can observe that, in the case of no noise, our proposed method captures the exact features of data for $\lambda = 1.5$, while the classical least squares method does not capture the data exactly. Moreover, when we increase the noise level, our proposed method is also accurate within two points significant digits.
\end{itemize}
\begin{table}[h!]
	\begin{center}
		\begin{tabular}{|c| c| l| c|c|}
			\hline
			& $\lambda$ & $\boldsymbol{y(x) = x^{1.5}}$ & Total Error $E^{D}$  & CPU time (in second)\\
			\hline
			\multirow{3}{*}{No Noise} & 1.50 & $y(x) \approx 0.0000+1.0000x^{1.5}$ & 8.20$e$-28 & 0.0149 \\
			\cline{2-5}
			& 1.25 & $y(x) \approx -11.1149+2.3609x^{1.25}$ & 2.02$e$-0 & 0.0145\\
			\cline{2-5}
			& 1.00 & $y(x) \approx -27.7817+8.1671x$ & 8.17$e$-0 & 0.0143\\
			\hline
			\multirow{3}{*}{$5 \%$ Noise} & 1.50 & $y(x) \approx -0.0001+1.0000x^{1.5}$ & 2.20$e$-2 & 0.0135 \\
			\cline{2-5}
			& 1.25 & $y(x) \approx -11.1136+2.3609x^{1.25}$ & 2.05$e$-0  & 0.0212\\
			\cline{2-5}
			& 1.00 & $y(x) \approx -27.7792+5.7898x$ & 8.19$e$-0 & 0.0141\\
			\hline
			\multirow{3}{*}{$10 \%$ Noise} & 1.50 & $y(x) \approx -0.0013+1.0000x^{1.5}$ & 8.80$e$-2 & 0.0142 \\
			\cline{2-5}
			& 1.25 & $y(x) \approx -11.1162+2.3610x^{1.25}$ & 2.10$e$-0 & 0.0142\\
			\cline{2-5}
			& 1.00 & $y(x) \approx -27.7859+5.7902x$ & 8.26$e$-0 & 0.0198\\
			\hline
		\end{tabular}
		\caption{Comparison of the least squares error $E^{D}$ in Example \ref{6ex2} corresponding to $n = 1$ and  with different value of $\lambda$.}\label{6tb2}
	\end{center}
\end{table}
\begin{ex}\label{6ex3}
	In this Example, we consider the data given in the Table \ref{6tb3}.
	\begin{table}[h!]
	\begin{center}
		\begin{tabular}{|c|c|c|c|c|c|}
			\hline 
			$x$ & $0.0000$ & $0.2500$ & $0.5000$ & $0.7500$ & $1.0000$ \\
			\hline
			$y$ with no noise & $0.0000$ & $0.1340$ & $0.3660$ & $0.6589$ & $0.6589$ \\
			\hline 
			$y$ with $5\%$ noise & $-0.0062$ & $0.0745$ & $0.0705$ & $0.0709$ & $0.0336$ \\
			\hline 
			$y$ with $10\%$ noise  & $-0.0062$ & $0.0745$ & $0.0705$ & $0.0709$ & $0.0336$ \\
			\hline 
		\end{tabular}
	\caption{Data with no noise, $5\%$ noise and $10\%$ noise.}\label{6tb3}
	\end{center}
\end{table}
\end{ex}
Here, data given in the Table \ref{6tb3} with no noise, $5\%$ noise and $10\%$ noise. For this Example, we search a best fit curve of the form of $P_{1}(x;\lambda) = a_{0}+a_{1}x^{\lambda} \in M_{1}^{\lambda}$. After implementation of the proposed method discussed in Section \ref{6s4}, the results are described in detail as below:
\begin{itemize}
	\item Table \ref{6tb7} shows the least squares error $E^{D}$ by the proposed method for $n = 1$ and different value of $\lambda$.
	\item Also, Table \ref{6tb7} displays the results for the noisy data set given in the Table \ref{6tb3}. Table results show the our proposed method is robust to noise.
	\item From Table \ref{6tb7}, one can observe that, in the case of no noise, we get more accurate results for $\lambda = 1.5$, while for noisy data, we get the precise result for a different choice of $\lambda$.
\end{itemize}
\begin{table}[h!]
	\begin{center}
		\begin{tabular}{|c| c| l| c|c|}
			\hline
			& $\lambda$ & $P_{1}(x;\lambda)$ & Total Error $E^{D}$ & CPU time (in second)  \\
			\hline
			\multirow{3}{*}{No Noise} & 1.50 & $ 0.0071+0.9977x^{1.5}$ & 1.3042$e$-4 & 0.0123\\
			\cline{2-5}
			& 1.00 & $-0.0732+0.0164x$ & 1.6400$e$-2 & 0.0109 \\
			\cline{2-5}
			& 0.5 & $-0.1409+0.9318x^{0.5}$ & 1.2320$e$-1 & 0.0115\\
			\hline
			\multirow{3}{*}{$5 \%$ Noise} & 1.50 & $0.0420+0.0157x^{1.5}$ & 4.7000$e$-3 & 0.0115 \\
			\cline{2-5}
			& 1.00 & $0.0335+0.0304x$ & 4.3000$e$-3 & 0.0107\\
			\cline{2-5}
			& 0.50 & $ 0.0159+0.0533x^{0.5}$ & 3.1000$e$-3 & 0.0131\\
			\hline
			\multirow{3}{*}{$10 \%$ Noise} & 1.50 & $ 0.1789-0.2227x^{1.5}$ & 3.5000$e$-2 & 0.0111 \\
			\cline{2-5}
			& 1.00 & $ 0.2139-0.2596x$ & 2.5400$e$-2 & 0.0113\\
			\cline{2-5}
			& 0.50 & $0.2705-0.3033x^{1.5}$ & 1.1130$e$-2 & 0.0111\\
			\hline
		\end{tabular}
		\caption{Comparison of the least squares error $E^{D}$ in Example \ref{6ex3} corresponding to $n = 1$ and  with different value of $\lambda$.}\label{6tb7}
	\end{center}
\end{table}
\begin{ex}\label{6ex7}
	In this Example, we consider the function $y(x) = x^{0.75}$ defined on the interval $[0,1]$ to generate the data set of length $30$.
\end{ex}
To demonstrate the advantage of choosing the weight functions, we consider Example \ref{6ex7}. For this, we divide the interval $[0, 1]$ into two parts $[0, 0.8]$ and $[0.8, 1]$. In $[0, 0.8]$, we take the data sets as it is. However, in $[0.8, 1]$, we add $10\%$ Gaussian noise in the data set. We generate fractional orthogonal polynomials concerning the two different weight functions. The outcome of these numerical experiments are described below:
\begin{itemize}
	\item In the first case, if one can generate the fractional orthogonal polynomials up to order $1$ with respect to weight function $W(x) = (1-x)$ using Theorem \ref{64thm1}, they get $L_{0}(x;0.75,W) = 1$ and $L_{1}(x;0.75,W) = x^{0.75}-0.5020$.
	\item In the second case, if one can generate the fractional orthogonal polynomials up to order $1$ with respect to weight function $W(x) = 1$, they get $L_{0}(x;0.75,W) = 1$ and $L_{1}(x;0.75,W) = x^{0.75}-0.7870$.
	\item The least-squares error $E^{D}$ in the first case is 2.020$e$-2, while in the second case, we get 4.2040$e$-1.
	\item From the above discussion one can conclude that choosing the appropriate weight function according to the data in the least squares method provides better results.
\end{itemize}
\begin{ex}\label{6ex9}
	In this Example, we will implement our proposed method to predict the value of American put options, where the risk-neutral stock price process satisfies the following stochastic differential equation:
	\begin{equation}\label{6e1}
		dS(x) = rS(x)dx+\sigma S(x)W(x),
	\end{equation}
where $r$ and $\sigma$ are constant, and $W(x)$ is the standard Brownian process. Here the variable $x$ is denoted the time.
\end{ex}
The well known solution of the Equation \eqref{6e1} is
\begin{equation}\label{6e2}
	S(x) = S_{0}\exp ((r-\frac{1}{2}\sigma^{2})x+\sigma W(x)),
\end{equation}
where $S_{0}$ is the initial stock price. The least square regression analysis is an essential part of machine learning. Therefore, we are interested in the predicted value of the American put option using our proposed method. We assume $S_{0} = 38,~r = 5\%,~\sigma = 71\%$, options strike price is $48$ and possible exercise time is $60$ days (see \cite{huang2009least} for details). Also, in the \cite{huang2009least}, the authors give the value of the American put option for the same data described above, which is $10.822$. After implementation of our proposed method, we get the following results:
\begin{itemize}
	\item Firstly we generate the data of length $60$ from the Equation \eqref{6e2} in the interval $[0,\frac{1}{6}]$ and then we fit the data using modified least square method in the space $M_{2}^{\lambda}$. The data is stochastic in nature, therefore, we use $10000$ simulation for prediction.
	\item Using the algorithm described in \cite{huang2009least} combined with the modified least square method, the value of predicated American put options are shown in Table \ref{6tb9} for different value of $\lambda$.
	\item From Table \ref{6tb9}, one can easily observe that the predicted value of American put options in the case of $\lambda = 0.75$ is much closer to the given value of the put options than other values of $\lambda$. 
\end{itemize}
\begin{table}[h!]
	\begin{center}
		\caption{Prediction of American put options in Example \ref{6ex9} with different value of $\lambda$.}
		\label{6tb9}
		\begin{tabular}{ |c|c|c|c|c| }
			\hline
			$\lambda$ & 0.25 & 0.5 & 0.75 & 1.00  \\
			\hline
			Predicated Value & 10.743 & 10.730 & 10.790  & 10.714 \\
			\hline 
			CPU Time (in seconds) & 10.720 & 10.723 & 10.756 & 10.856\\
			\hline 
		\end{tabular}
	\end{center}
\end{table}
\section{Application of modified least squares method}\label{6s6}
The modified least squares method will have many practical applications in physics, finance, and other engineering problems. In this Section, we have demonstrated the application of the modified least squares method in particular areas like solving fractional differential/integral equations and in machine learning.
\subsection{Application in solving fractional differential/integral equations}
\begin{figure}[h!]
	\begin{center}
		\includegraphics[scale=0.75]{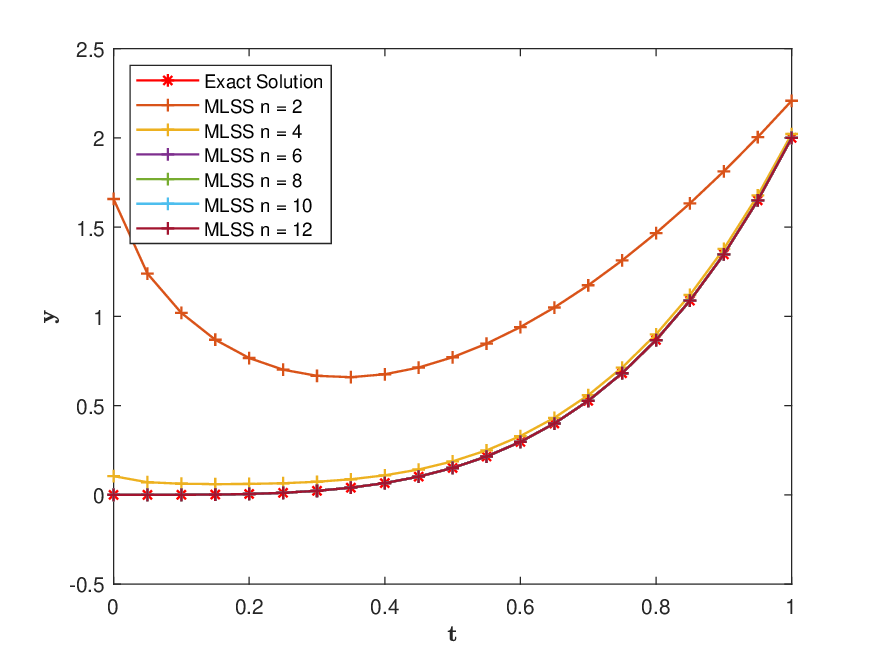}
	\end{center}
	\caption{Exact solution and approximate solution for Example \ref{6ex8} with $\lambda = 0.75$ and different values of $n$.}
	\label{6fig1}
\end{figure}

The theory of non-integer derivatives is an emerging topic of applied mathematics, which attracted many researchers from various disciplines. The non-local properties of fractional operators attract a significant level of intrigue in the area of fractional calculus. It can give an excellent way to deal with complex phenomena in nature, such as biological systems, control theory, finance, signal and image processing, sub-diffusion and super-diffusion process, viscoelastic fluid, electrochemical process, and so on (see \cite{herrmann2013folded,odzijewicz2013noether,singh2021wavelet,mehandiratta2019existence} and references therein). The main advantage of fractional differential/integral equations is that it provides a powerful tool for depicting the system with memory, long-range interactions and hereditary properties of several materials instead of the classical differential/integral equations in which such effects are difficult to incorporate. The fractional differential equations are equivalent to the Volterra's second kind integral equations for the specific choice of kernel. Consider the following Volterra second kind integral equation of the form
\begin{equation}\label{6eq1}
	y(x) = a + \int_{0}^{x}k(x,t)f(t,y(t))dt,
\end{equation}
where $a \in \mathbb{R},~f(x,t)$ to be a continuous function whereas $k(x,t)$ may be singular. When $k(x,t) = (x-t)^{\alpha-1},~0 < \alpha <1$, then Equation \eqref{6eq1} is equivalent to the following fractional differential equation
\begin{equation}
	\begin{aligned}
		\textsubscript{C}D_{0,x}^{\alpha}y(x) &= \Gamma(\alpha)f(x,y(x)), \\
		y(0) &= a.
	\end{aligned}
\end{equation}
However, in many cases, it is not possible to find the exact solution for fractional differential equations. Therefore, it is essential to acquire its approximate solution by using some numerical methods. In the literature there are a paper for solving fractional differential equations using least squares method based on polynomials $P_{n}(x;1) \in M_{n}^{1}$ \cite{rahimkhani2017fractional}. But the fractional derivative of the polynomials $P_{n}(x;1)$ does not belong to the space $M_{n}^{1}$. Therefore, we introduce some extra error while solving the fractional differential equations using least squares method based on $P_{n}(x;1)$. For example, $n=2$, consider the space 
\begin{equation*}
M_{2}^{\lambda} := \{1,x^{\lambda},x^{2\lambda}\}.
\end{equation*}
The fractional derivative (in Caputo sense) of order $\alpha$ of any polynomial $P_{i}(x;\lambda) \in M_{2}^{\lambda},i=0,1,2$ are also in the space $M_{2}^{\lambda}$ for the choice $\alpha = \lambda$. Hence, we feel that when solving the fractional differential using the least squares method in the space $M_{n}^{\lambda}$ is beneficial than the space $M_{n}^{1}$. 
\begin{ex}\label{6ex4}
	Consider the following fractional differential/integral equation
	\begin{equation}
	\textsubscript{C}D_{0,x}^{\alpha}y(x) = \frac{1}{\Gamma(2-\alpha)}x^{1-\alpha},~~0 < x \leq 1,
	\end{equation}
	with the analytical solution $y(x) = x$ when $y(0) = 0$ and $f(x) = \frac{1}{\Gamma(2-\alpha)}x^{1-\alpha}$. 
\end{ex}
We are using the modified least squares method to solve Example \ref{6ex4} in space $M_{2}^{\lambda}$. Let the solution of fractional differential equation $y(x)$ be approximated by the polynomial $a_{0}+a_{1}x^{\lambda}+a_{2}x^{2\lambda}$. In this case, we define the residual error as $R(x,a_{0},a_{1},a_{2};\lambda) = f(x)-\textsubscript{C}D_{0,x}^{\alpha}(a_{0}+a_{1}x^{\lambda}+a_{2}x^{2\lambda})$. So, our error function becomes
\begin{equation}
\begin{aligned}
E^{C}(a_{0},a_{1},a_{2};\lambda) &= \int_{0}^{1}(R(x,a_{0},a_{1},a_{2};\lambda))^{2}dx.
\end{aligned}
\end{equation} 
For fix $n = 2$, the least square error $E^{C}$ for Example \ref{6ex4} has been shown in Table \ref{6tb4} for $\alpha = 0.5$ and various values of $\lambda$. From Table \ref{6tb4}, one can observe that when $\lambda = \alpha$, we get the best result. This will happen because fractional derivatives of $x^{\lambda}$ and $x^{2\lambda}$ belong to the space $M_{2}^{\lambda}$ when $\lambda = \alpha$.
\begin{table}[h!]
	\begin{center}
		\caption{Comparison of the least square error $E^{C}$ in Example \ref{6ex4} corresponding to $\alpha = 0.5$ and  with different value of $\lambda$.}
		\label{6tb4}
		\begin{tabular}{ |c|c|c|c|c|c| }
			\hline
			$\lambda$ & 0.5 & 0.75 & 1.00 & 1.25 & 1.50 \\
			\hline
			$E^{C}$ & 0 & 6.11$e$-4 & 5.19$e$-4 & 2.70$e$-3 & 8.60$e$-3  \\
			\hline 
		\end{tabular}
	\end{center}
\end{table}
\begin{ex}\label{6ex8}
	In this Example, we consider the following multi-term fractional/integral differential equation
	\begin{equation}
		\textsubscript{C}D_{0,x}^{\alpha_{1}}y(x)+\textsubscript{C}D_{0,x}^{\alpha_{1}}y(x)+y(x) = f(x),~~ 0 < x \leq 1,
	\end{equation}
	with the exact solution $y(x) = x^{3.5}+x^{4}$ when $y(0) = 0$ and $f(x) = x^{4}+x^{3.5}+\frac{\Gamma(5)}{\Gamma(5-\alpha_{1})}x^{4-\alpha_{1}}+\frac{\Gamma(4.5)}{\Gamma(4.5-\alpha_{1})}x^{3.5-\alpha_{1}}+\frac{\Gamma(5)}{\Gamma(5-\alpha_{2})}x^{4-\alpha_{2}}+\frac{\Gamma(4.5)}{\Gamma(4.5-\alpha_{2})}x^{3.5-\alpha_{2}},$ where $\alpha_{1} = 0.5$ and $\alpha_{2} = 0.25$.
\end{ex}
We are using the modified least squares method to solve Example \ref{6ex8} in space $M_{n}^{\lambda}$. Let the solution of fractional differential equation $y(x)$ be approximated by the M\"untz-Legendre polynomials $\sum_{i=0}^{n}a_{i}L_{i}(x;\lambda)$. In this case, we define the residual error as 
\begin{equation*}
	\begin{aligned}
		R(x,a_{0},a_{1},\ldots,a_{n};\lambda)& = \textsubscript{C}D_{0,x}^{\alpha_{1}}(\sum_{i=0}^{n}a_{i}L_{i}(x;\lambda))+\textsubscript{C}D_{0,x}^{\alpha_{2}}(\sum_{i=0}^{n}a_{i}L_{i}(x;\lambda))\\
		&	+\sum_{i=0}^{n}a_{i}L_{i}(x;\lambda)-f(x)  +\sum_{i=0}^{n}a_{i}L_{i}(0;\lambda)-y(0).
	\end{aligned}
\end{equation*} So, our error function becomes
\begin{equation}
	\begin{aligned}
		E^{C}(a_{0},a_{1},\ldots,a_{n};\lambda) &= \int_{0}^{1}(R(x,a_{0},a_{1},\ldots,a_{n};\lambda))^{2}dx.
	\end{aligned}
\end{equation} 
The least square error $E^{C}$ and absolute error (A.E.) at $x = 1.00$ for Example \ref{6ex8} and corresponding CPU time have been shown in Table \ref{6tb8} for various values of $\lambda$ and $n$. One can easily observe that for $\lambda = 1$, the approximate solution does not capture the exact features involved in the solutions of the Example \ref{6ex8} for any value of $n$. However, for some values of $\lambda \neq 1$, the approximate solution captures the exact features involved in the solutions of the Example \ref{6ex8}. Thus, for the conclusion of the results of this Example, we can say that one can get a good approximation for the solution of fractional differential equations in the case of $ \lambda \neq 1 $. Figure \ref{6fig1} shows the trajectory of the exact solution and modified least square solution (MLSS) for Example \ref{6ex8} with $\lambda = 0.75$ and different values of $n$. Figure \ref{6fig1} shows that our approximate solution converges to the exact solution when we increase the value of $n$.
\begin{table}[h!]
\begin{center}
	\caption{Comparison of the least square error $E^{C}$ and A.E. in Example \ref{6ex8} corresponding for different values of $\lambda$ and $n$.}
	\label{6tb8}
	\begin{tabular}{|c|c|c|c|c|c|}
		\hline
		\multicolumn{2}{|c|}{$\lambda$} & $0.50$ & $0.75$ & $1.00$ & $1.25$\\
		\hline
		\multirow{2}{*}{$n=2$} &
		$E^c$ & 5.96$e$-1 & 3.19$e$-1  & 1.29$e$-1 & 4.31$e$-2 \\
		\cline{2-6}
		& A.E. & 2.22$e$-0 & 2.08$e$-1 & 3.28$e$-2 & 2.67$e$-4 \\
		\cline{2-6}
		& CPU Time (in seconds) & 2.69 & 3.08  & 3.28 & 3.16\\
		\hline
		\multirow{2}{*}{$n=4$} &
		$E^c$ & 7.93$e$-3 & 2.33$e$-4 & 5.10$e$-7  & 1.85$e$-9 \\
		\cline{2-6}
		& A.E. & 2.09$e$-0 & 2.19$e$-2 & 4.27$e$-4 & 3.45$e$-5 \\
		\cline{2-6}
		& CPU Time (in seconds) & 5.89 & 9.33 & 10.67 & 13.90\\
		\hline
		\multirow{2}{*}{$n=6$} &
		$E^c$ & 4.59$e$-6  & 2.57$e$-11 & 9.79$e$-10 & 3.05$e$-10\\
		\cline{2-6}
		& A.E. & 1.72$e$-1 & 8.81$e$-6 & 3.27$e$-5 & 1.18$e$-5\\
		\cline{2-6}
		& CPU Time (in seconds) & 16.55 & 31.06 & 34.29 & 37.73\\
		\hline
		\multirow{2}{*}{$n=8$} &
		$E^c$ & 3.08$e$-45 & 6.61$e$-14 & 1.53$e$-11 & 2.16$e$-11 \\
		\cline{2-6}
		& A.E. & 4.40$e$-16  & 1.34$e$-6 & 5.79$e$-6 & 3.98$e$-6 \\
		\cline{2-6}
		& CPU Time (in seconds) & 39.55 & 75.15 & 77.97 & 81.87 \\
		\hline
		\multirow{2}{*}{$n=10$} &
		$E^c$ & 2.69$e$-47 & 2.27$e$-15 & 6.45$e$-13 & 2.35$e$-12\\
		\cline{2-6}
		& A.E. & 4.40$e$-16 & 3.12$e$-7 & 1.53$e$-6 & 1.59$e$-6 \\
		\cline{2-6}
		& CPU Time (in seconds) & 80.23  & 149.29 & 183.23 & 241.30\\
		\hline
	\end{tabular}
\end{center}
\end{table}
\subsection{Application in machine learning}
Machine learning is a field of artificial intelligence which gives computer systems the ability to learn using available data. Recently machine learning algorithms become very popular for analyzing of data and make predictions. The least-squares method is widely used in machine learning to analyze data for regression analysis and classification. In particular, in the regression analysis, the goal of people is to plot a best fit curve or line between data. If someone is interested in discovering the best fit line for one variable using the least squares method for the data, then, in this case, the form of polynomials (hypothesis/model) is
\begin{equation*}
P_{1}(x;1) = a_{0}+a_{1}x.
\end{equation*}
We have already shown the advantage of the modified least squares method in Section \ref{6s5} while searching the best fit curve between data. Also, some data need to have features vector like $x^{\lambda}$, for example $x^{0.5}$. Therefore, we must choose the polynomials (hypothesis/model) of the form 
\begin{equation*}
P_{2}(x;0.5) = a_{0}+a_{1}x^{0.5}.
\end{equation*}
\begin{ex}\label{6ex5}
	Consider the data in Table \ref{6tb5}, which is related to pharmaceutical sales of some company.
 \begin{table}[h!]
 	\begin{center}
 		\caption{Pharmaceutical sales of some company \cite{abazid2018least}}
 		\label{6tb5}
 		\begin{tabular}{ |c|c|c|c|c|c| }
 			\hline
 			Years & 2014 & 2015 & 2016 & 2017 & 2018 \\
 			\hline
 			Sales & 10000 & 21000 & 50000 & 70000 & 71000 \\
 			\hline 
 		\end{tabular}
 	\end{center}
 \end{table}
\end{ex}
For Example \ref{6ex5}, we fit the data with the $P_{1}(x;\lambda) = a_{0}+a_{1}x^{\lambda} \in M_{1}^{\lambda}$. To simplify the calculation, the years in Table \ref{6tb5} are replaced by the coded values. For example, $2014$ is $1$, $2015$ is $2$, and so forth on. We use the data from $2014$ to $2017$ to fit the curve $a_{0}+a_{1}x^{\lambda}$ and predict the sales in $2018$. After implementation of the proposed method discussed in Section \ref{6s4}, the results are described in the Table \ref{6tb6}. From the Table \ref{6tb6}, one can observe that for $\lambda = 0.5$, the predicted sales in $2018$ is close to the real value in $2018$.
\begin{table}[h!]
	\begin{center}
		\caption{Prediction of pharmaceutical sales of some company in the year $2018$ with different values of $\lambda$}
		\label{6tb6}
		\begin{tabular}{ |c|c|c|c|c|c| }
			\hline
			$\lambda$ & 0.5 & 0.75 & 1.00 & 1.25 & 1.5 \\
			\hline
			Predicted sales in 2018 & 69692 & 80546 & 90000 & 98307 & 105870 \\
			\hline 
		\end{tabular}
	\end{center}
\end{table}
\begin{ex}\label{6ex10}
	In this Example, we will fit the data generated by the solution of the fractional differential equation, which describes the dynamics of the world population.
	\begin{equation}\label{6enew}
		\begin{aligned}
			\textsubscript{C}D_{0,x}^{\alpha}y(x) &= Py(x),~~0 < x \leq 1, \\
			y(0) &= y_{0} 
			\end{aligned}
	\end{equation} The solution of the fractional differential equations is given by $y(x) = y_{0}E_{\alpha}(Px^{\alpha})$, where $y_{0}$ is the population of the world at initial time, $E_{\alpha}(\cdot)$ is the Mittage-Leffler function, $P$ is the production rate, and $\alpha$ is the fractional order of the model. 
\end{ex}
In \cite{almeida2016modeling}, authors are finding the values of $P$ and $\alpha$ using the world population data from $1910$ to $2010$, which are $1.3502e-2$ and $1.39$, respectively. For this Example, we generating the data of $y(x)$ with the help of the exact solution of the fractional differential equation \eqref{6enew} at $11$ equispaced points on the interval $[0,1]$ for the values of $P$ and $\alpha$. Here, we consider the two cases: first, we approximate $y(x)$ with the  $\sum_{i=0}^{n}a_{i}x^{i\lambda}$, while in the second case, we approximate $y(x)$ with the $\sum_{i=0}^{n}b_{i}P_{i}(x;\lambda)$, where $P_{i}(x;\lambda),~i = 0,1,\ldots,n,$ are the orthogonal fractional polynomials with respect to data $\{x_{j},~j=0,1,\ldots,10\}$, and we generated $P_{i}(x;\lambda),~i=0,1,\ldots,n,$ with help of the Theorem \ref{64thm2}. In the first case for finding the value of parameter $a_{i},~i=0,1,\ldots,n$, we have to solve the system of equations because $\{1,x^{\lambda},\ldots,x^{n\lambda}\}$ are non-orthogonal fractional polynomials with respect to data $\{x_{j},~j=0,1,\ldots,10\}$, and for a large value of $n$, we may end with the ill-conditioned coefficient matrix. In the second case, we are directly finding the parameter values $b_{i},~i=0,1,\ldots,n$ using Equation \eqref{63e10}. Table \ref{6tb10} demonstrates the absolute error (A.E.) at $x = 0.55$ with different values of $\lambda$ and $n$ for both cases. From Table \ref{6tb10}, one can observe that in the case of $\lambda = 1.39$, we get the better results compared to other values of $\lambda$ for each $n$. This happens because the exact data is generated for $\alpha = 1.39$.
\begin{table}[h!]
	\begin{center}
		\caption{A.E. at $x = 0.55$ for Example \ref{6ex10} with different values of $\lambda$ and $n$.}
		\label{6tb10}
		\begin{tabular}{ |c|c|c|c|c| }
			\hline
			\multicolumn{5}{|c|}{For non-orthogonal fractional polynomials} \\
			\hline 
			$n$ & $\lambda = 0.50$ & $\lambda = 1.00$ & $\lambda = 1.50$ & $\lambda = 1.39$  \\
			\hline
			2 &  1.96$e$-4 & 4.16$e$-5 & 4.78$e$-7 & 6.68$e$-10 \\
			\hline 
			3 &  8.79$e$-7 & 1.61$e$-5 & 1.50$e$-5 & 6.86$e$-13 \\
			\hline 
			4 &  4.19$e$-7 & 7.01$e$-6 & 1.62$e$-6 & 5.10$e$-15 \\
			\hline
			5 &  2.64$e$-8 & 2.00$e$-6 & 4.79$e$-6 & 5.11$e$-15 \\
			\hline  
			6 &  3.61$e$-9 & 1.35$e$-6 & 8.31$e$-7 & 2.21$e$-13 \\
			\hline 
			\multicolumn{5}{|c|}{For orthogonal fractional polynomials} \\
			\hline 
			$n$ & $\lambda = 0.50$ & $\lambda = 1.00$ & $\lambda = 1.50$ & $\lambda = 1.39$  \\
			\hline
			2 &  4.61$e$-4 & 2.01$e$-3 & 1.08$e$-4 & 8.36$e$-9 \\
			\hline 
			3 &  1.05$e$-5 & 9.55$e$-5 & 8.68$e$-5 & 4.76$e$-12 \\
			\hline 
			4 &  1.98$e$-6 & 6.29$e$-5 & 1.00$e$-4 & 3.18$e$-15 \\
			\hline
			5 &  2.63$e$-7 & 4.93$e$-5 & 1.54$e$-4 & 4.83$e$-15 \\
			\hline  
			6 &  5.57$e$-8 & 4.24$e$-5 & 3.12$e$-4 & 2.09$e$-16 \\
			\hline 
		\end{tabular}
	\end{center}
\end{table}
\section{Conclusions and future work}\label{6s7}
The introductory Section shows the demands of the least squares method in various fields. So the modification in the least squares method is the demand of time due to its application. The main idea of this work has been to modify the least squares method using the space $M_{n}^{\lambda}$. The numerical results for test Examples have been reported to show the efficiency of the modified least squares method over the classical least squares method. We can use the current work in the support vector machines in the future.
Support vector machines are part of machine learning to analyze data for classification and regression analysis. In most cases, data are non-linear. So, we find some non-linear transformation $\phi$ that can be mapped the data onto high-dimensional feature space. The transformation is chosen in such a way that their dot product leads the kernel style function
\begin{equation*}
K(x,x_{i}) = \phi(x).\phi(x_{i}).
\end{equation*}
If we choose the polynomial classifiers \cite{graf1997polynomial} of degree $2$  and we have $n$ data set $x_{i},~i=1,\ldots,n$, in this case
\begin{equation*}
\phi(x) = (1 ~x_{1}~ x_{2} \ldots x_{n}~ x_{1}^{2}~ x_{1}x_{2} \ldots x_{2}^{2} \ldots x_{n}^{2})^{T}.
\end{equation*}
However, one can use fractional polynomial classifiers instead of classical polynomial classifiers for more accurate results. Further, one can also use fractional orthogonal polynomials as an activation function in the neural networks to avoid the vanishing gradient problem.
\section*{Acknowledgements}
The First author acknowledges the support provided by University Grants Commission (UGC), India, under the grant number $20/12/2015$(ii)EU-V. The second author acknowledges the support provided by the SERB India, under the grant number SERB/F/$3060/2021-2022$. The third author acknowledges the financial support from the North-Caucasus Center for Mathematical Research under agreement number $075-02-2021-1749$ with the Ministry of Science and Higher Education of the Russian Federation.
\section*{Declaration of Competing Interest}
The authors declare that they have no known competing financial interests or personal relationships that could have appeared to influence the work reported in this paper.
\bibliographystyle{plain}
\bibliography{main}
\end{document}